\documentclass[12pt]{amsart}
\usepackage{pgfplots}
\usepackage{float}
\usepackage[utf8]{inputenc}
\usepackage{graphicx}
\usepackage{color}
\usepackage{mdwlist}
\usepackage{amssymb}
\usepackage{subfigure}
\usepackage{relsize}
\usepackage{graphics}
\allowdisplaybreaks

\usepackage{cleveref}

\newcommand{\beq}{\begin{eqnarray*}}
\newcommand{\feq}{\end{eqnarray*}}
\newcommand{\beqn}{\begin{eqnarray}}
\newcommand{\feqn}{\end{eqnarray}}

\newcommand{\RN}[1]{%
  \textup{\uppercase\expandafter{\romannumeral#1}}%
}

\newtheorem{theorem}{Theorem}[section]
\newtheorem*{theorem*}{Theorem}

\theoremstyle{definition}

\theoremstyle{remark}

\numberwithin{equation}{section}
\begin{document}
\title[]
{On critical thresholds for hyperbolic balance law systems}% with a constant background}

\author{Manas Bhatnagar and Hailiang Liu}
\address{Department of Mathematics and Statistics, University of Massachusetts Amherst, Amherst, Massachusetts 01003}
\email{mbhatnagar@umass.edu}
\address{Department of Mathematics, Iowa State University, Ames, Iowa 50010}
\email{hliu@iastate.edu} 
%\address{Department of Mathematics, University of South Carolina, Columbia, South Carolina 29208}
% Address of record for the research reported here
%\address{Department of Mathematics, Iowa State University, Ames, Iowa 50010}
%\email{tan@math.sc.edu}
 %35L65 : Conservation laws,  35L67 : Shocks and singularities
%35B30  Dependence of solutions on initial and boundary data, parameters
\begin{abstract}  We review the theoretical development in the study of critical thresholds for hyperbolic balance laws. The emphasis is on two classes of systems: Euler-Poisson-alignment (EPA) systems and hyperbolic relaxation systems. We start with an introduction to the `Critical Threshold Phenomena' and study some nonlocal PDE systems, which are important from modeling point of view. 
\end{abstract}
\maketitle

\section{Introduction}
\label{sec:1}
\iffalse
Use the template \emph{chapter.tex} together with the document class SVMono (monograph-type books) or SVMult (edited books) to style the various elements of your chapter content.

Instead of simply listing headings of different levels we recommend to let every heading be followed by at least a short passage of text. Further on please use the \LaTeX\ automatism for all your cross-references and citations. And please note that the first line of text that follows a heading is not indented, whereas the first lines of all subsequent paragraphs are.
\fi
For first order hyperbolic conservation laws, it is generic that the solutions lose smoothness even if the initial data is smooth, \cite{Da16, Lax64}. However, addition of source terms can balance this `breaking' and result in global-in-time smooth solutions for a large class of initial data. 
%There has been an enormous amount of papers related to the study of global behavior of solutions to hyperbolic balance laws. 
For the question of global behavior of strong solutions, the choice of the initial data and/or damping forces is decisive. The classical stability analysis can fail either for large perturbations of some wave patterns or when the steady state solution may be only conditionally stable due to the weak dissipation in the system, see for example \cite{GHZ17,Y01}. On the other hand, the notion of critical threshold (CT) has been shown to be powerful in describing the conditional stability for underlying physical problems, and the associated phenomena does reflect the delicate balance among various forcing mechanisms, \cite{BL20,BL201,CCTT16,ELT01}. 

An example to illustrate this is the pressureless Euler-Poisson  (EP) system that consists of the continuity equation, Burgers' equation with an electric source through a potential, and the Poisson equation for the potential,
\begin{align}
\label{pressurelesseuler}
\begin{aligned}
    & \rho_t + (\rho u)_x = 0,\\
    & u_t + uu_x = -\phi_x,\\
    & -\phi_{xx} = \rho,
\end{aligned}
\end{align}
with smooth initial data $(\rho_0\geq 0, u_0)\in C^1(\mathbb{R})\times C^1(\mathbb{R})$.
Taking the spatial derivative of the second equation and setting $g(t,X):=u_x(t,X)$ for $\frac{dX}{dt} = u(t,X)$, we can obtain an ODE system along the characteristics,
\begin{align*}
&\frac{d (\rho(t,X))}{dt} = -\rho g,\\
& \frac{d (g(t,X))}{dt} = -g^2 + \rho.
\end{align*}
We have omitted the parameter for the ODE, that is a consequence of the method of characteristics, to avoid excess notation. For global well-posedness of \eqref{pressurelesseuler}, the issue now has reduced to ensuring the all-time-boundedness of $\rho,g$ as solutions of the above ODE system. If we assume the initial density to be nowhere zero, then the first ODE implies that the density remains positive for all times. Using simple calculations, one has,
\begin{align*}
& \frac{d}{dt}(g/\rho) = 1,\qquad \frac{d}{dt}(1/\rho) = \frac{g}{\rho}.
\end{align*}
Hence, 
$$
\frac{1}{\rho(t,X(t))} = \frac{t^2}{2} + t\frac{g(0)}{\rho(0)} + \frac{1}{\rho(0)}.
$$
From this expression, we can conclude that $\rho(t)$ is finite for all time if and only if $g(0)\geq -\sqrt{2\rho(0)}$. 
%Moreover, if $g(0)<-\sqrt{2\rho(0)}$, then $\rho(t,X(t))\to \infty$ in finite time. 
Therefore, for \eqref{pressurelesseuler}, there is global classical solution if and only if for all $x\in\mathbb{R}$, we have $(u_{0x}(x),\rho_0(x))\in \{(g,\rho):g\geq-\sqrt{2\rho},\rho\geq 0\}$. The question of global-in-time existence vs finite-time-breakdown has now been reduced to whether the initial data crosses a certain threshold, which here, is the curve,  $\{(g,\rho):g=-\sqrt{2\rho},\rho\geq 0\}$ (on the phase space of $(u_{0x},\rho_0)$). This is the \textit{CT curve}. It divides the initial data into two sets: the \textit{subcritical} region (the good part) and the \textit{supercritical} region (the bad part).
%Moreover, if there exists some $x_\ast$ such that $u_{0x}(x_\ast)<-\sqrt{\rho_0(x_\ast)}$, then 
%$$
%\lim_{t\to t_c^-}u_x(t,x_c) = -\infty = -\lim_{t\to t_c^-}\rho(t,x_c),
%$$ 
%for some $x_c\in\mathbb{R}, t_c>0$. 
%Consequently, this allows for an initial profile of velocity which can have negative slopes, resulting in a nontrivial set of initial profiles that results in global solutions. This is not the case for inviscid Burgers', where there can be no such initial data.

%For \eqref{pressurelesseuler}, we call the curve $\{(x,y):x=-\sqrt{2y},y\geq 0\}$ (on the phase space of $(u_{0x},\rho_0)$), the \textit{CT curve}. The set $\{(x,y):x=-\sqrt{2y},y\geq 0\}$ is the \textit{subcritical} region and the set $\{(x,y):x< x=-\sqrt{2y},y\geq 0\}$ is the \textit{supercritical} region. Here the threshold is a curve. However, in general, it could be a manifold. The threshold analysis as well as the manifolds have different levels of complexity based on the PDE system in question. 

The stimulating interplay of old and new ideas continue to lead to increased understanding as well as an ever-larger set of techniques with CT analysis for a larger array of problems.   The results on CT phenomena reviewed here are obtained in a series of papers \cite{BL211,BL212,BLT23,BL22}, which provide an account on recent developments of CT theory for a class of hyperbolic balance laws. The corresponding methods in obtaining these results will be briefly illustrated.

%In this chapter we review some of our results/techniques for Euler-like systems with regards to critical thresholds. 

\section{Euler-Poisson alignment dynamics}
\label{sec:2}
% Always give a unique label
% and use \ref{<label>} for cross-references
% and \cite{<label>} for bibliographic references
% use \sectionmark{}
% to alter or adjust the section heading in the running head
%Critical threshold analysis for equations with nonlocal terms is relatively difficult. Moreover, it is quite hard to obtain precise threshold curves/manifolds in such situations. One such system is 
The Euler-Poisson-alignment (EPA) system models phenomena where matter can be regarded as consisting of a continuum of moving particles, such as flow of charge or flocking of birds, see for example \cite{ELT01,TT14}. It consists of the following system of equations,
\begin{subequations}
\label{epamain}
\begin{align}
& \rho_t + (\rho u)_x = 0, \label{epamass}\\
& u_t + uu_x = -k\phi_x +\int_{\mathbb{R}}\psi(x-y)(u(t,y)-u(t,x))\rho(t,y)dy, \label{epamom}\\
& -\phi_{xx} = \rho -c, \label{epapoisson}
\end{align}
\end{subequations}
where  $x\in\mathbb{T}$ (periodic torus) or $x\in\mathbb{R}$, and $t>0$. The initial data is $\rho_0(x),u_0(x)$. Here, $k$ is the forcing coefficient, the sign of which signifies the type of forcing between particles. $c$ is the background term. The equation \eqref{epamom} has two types of forces: the alignment force modeled through the influence function $\psi$, and the electric force modeled through the potential $\phi$. The alignment force affects the collective motion of particles as to how they react to other particles around them. $\psi$ is \textbf{nonnegative} and \textbf{symmetric}. It models the strength and extent of the interaction. %The electric potential $\phi$ follows the Poisson equation, \eqref{epapoisson}. 
\iffalse
In the next few paragraphs, we give a short physical motivation for the potential $\phi$. Suppose the spatial domain in \eqref{epamain} to be all of $\mathbb{R}$ and $c=0$. Integrate \eqref{epapoisson} as follows,
$$
-\phi_x(t,x) = \int_{-\infty}^x \rho (t,y)\, dy.
$$
We assume the boundary condition for the Cauchy problem that $\lim_{y\to\pm\infty} \phi_x(t,y) = 0$. 
%It is worthe above expression implies the electric force in \eqref{epamom} is nonlocal. 
Similarly, integrating from $x$ to $\infty$ instead of $-\infty$ to $x$, we obtain,
$$
\phi_x(t,x) = \int_x^{\infty}\rho (t,y)\, dy.
$$
Taking difference, 
\begin{align*}
%\label{electricf}
-\phi_x = \int_{-\infty}^x \rho(t,y)\, dy  - \int_x^\infty \rho(t,y)\, dy.
\end{align*}
%where $M_0:=\int_{\mathbb{R}}\rho_0 dy$, which is conserved for all time owing to \eqref{epamass}. 
Plugging in the above expression into \eqref{epamom} and neglecting the alignment force ($\psi\ast(\rho u)-u\psi\ast\rho$) for the moment, we obtain,
$$
u_t + uu_x = k\left(\int_{-\infty}^x \rho(t,y)\, dy - \int_x^{\infty} \rho(t,y)\, dy\right) .
$$
\fi

$\phi$ is the potential that models electric or gravitational forces as and according when $k>0$ or $k<0$ respectively. 
%the system models gravitational forces and there is attraction instead of repulsion. 
The system \eqref{epamain}, with $\psi\equiv 0$ is called the Euler-Poisson (EP) system (illustrated in Section \ref{sec:1}) and has been a topic of intensive study by various researchers due to its ability to model numerous physical phenomena, \cite{Ja62,Ma86,MP90}.

When $k=0$ in \eqref{epamain}, it is called the Euler-alignment (EA) system. 
%Such a system was derived as a mean-field limit of an $N$-particle system in the limit as $N\to\infty$, see \cite{HaTa08} for a derivation. 
%EA systems are used to model various biological phenomena ranging from change of opinions to flocking of birds, \cite{CCTT16,TT14}. 
EA systems have been analyzed with different types of $\psi$. In \cite{Tan20}, the author obtained global existence results for $\psi\in L^1$. In \cite{BL22}, we relaxed the hypothesis and used a different, more elementary technique to arrive at the result. The gist of the result is mentioned in Theorem \ref{eal1}.

%Theorems \ref{gs1} and \ref{l1gs}  pertain to the full EPA system, \eqref{epamain} with $k>0$. Results for $k<0$ were derived by the authors in \cite{BL201}. The case when $k>0$ is more complicated than $k<0$ due to the presence of oscillatory solutions which are absent in the latter. The EP and EA systems have individually been studied for thresholds, \cite{BL20,BL201,CCTT16,CCZ16,ELT01,Tan20,TT14}. However, a full system like \eqref{epamain} has been difficult to analyze and there were no critical threshold results pertaining to this system before the authors' works. 

EPA systems with background ($c>0$) have to be studied on a periodic spatial domain and not on $\mathbb{R}$. This is owing to assumptions required for local well-posedness.  Therefore, we let the spatial variable space be $\mathbb{T}=[-1/2,1/2)$, the periodic torus. Local existence requires the assumption: $\int_{-1/2}^{1/2} \rho(t,y)-c\, dy = 0$.
This equality holds for all time if it holds initially. This is because mass is conserved by \eqref{epamass}. In view of this, we set, 
$$
c=\int_{-\frac{1}{2}}^{\frac{1}{2}}\rho_0(y)\, dy.
$$
Our first set of results is for $\psi\in L^\infty (\mathbb{T})$ having,
\begin{align}
\label{eq:psibounded}
& 0\leq \psi_{min}\leq\psi\leq \psi_{max}.
\end{align}
With the above assumption, \eqref{epamain} can be reformulated into a simpler system. Setting $G:= u_x + \psi\ast\rho$ with $G_0(x) := u_{0x}(x)+(\psi\ast\rho_0)(x)$, we obtain the following.
\begin{subequations}
\label{localexeqsys}
\begin{align}
& G_t + (Gu)_x = k(\rho -c), \label{localexeqsys2}\\
& \rho_t + (\rho u)_x = 0, \label{localexeqsys1}
%& u_x = G - \psi\ast\rho, \label{localexeqsys3}
\end{align}
\end{subequations}
with initial data $(G_0,\rho_0)\in H^s(\mathbb{T})\times (H^s(\mathbb{T})\cap L^1_+(\mathbb{T})) $, for $s>1/2$. The local existence for such system is known, \cite{CCTT16}. In particular, if initial data is smooth, then a smooth solution exists for some finite time. 
%Following term quantifies the strength of the electric force in the system and is helpful in stating our main results.
%\begin{align}
%\label{lambdaexp}
%& \lambda := 2\sqrt{\frac{k}{c}}.
%\end{align}
\begin{theorem}[Bounded alignment force]
\label{gs1}
Consider \eqref{localexeqsys} with repulsive electric force $k>0$ and
bounded alignment influence $\psi$ satisfying \eqref{eq:psibounded}. Set $\lambda := 2\sqrt{\frac{k}{c}}$.
Suppose the initial data $(G_0, \rho_0)$ is smooth and lies in the space mentioned above.
Then there exist sets $\Sigma_1,\Sigma_2,\Sigma_3$ such that,
\begin{enumerate}
    \item 
    \emph{Weak alignment} ($\psi_{max}<\lambda$):
    under the admissible condition 
    \begin{equation}\label{eq:weakacc}
    \psi_{max} - \psi_{min} < \frac{e^{\frac{\tan^{-1}\hat z}{\hat z}}\left(1 - e^{ -\frac{\pi}{\tilde z} - \frac{\pi}{\hat z} }\right) }{2\left(1 + e^{-\frac{\pi}{\tilde z}}\right)}\lambda,
    \end{equation}
    if the initial data lie in the subcritical region $\Sigma_1$, namely
    $$
    \big(G_0(x),\rho_0(x)\big)\in\Sigma_1, \quad \forall\,x\in\mathbb{T},
    $$
    then \eqref{localexeqsys} admits global-in-time classical solutions.
    \item
    \emph{Strong alignment} ($\psi_{min}\ge\lambda$):
    if the initial data lie in the subcritical region $\Sigma_2$, namely
    $$
   \big(G_0(x),\rho_0(x)\big)\in\Sigma_2, \quad \forall\,x\in\mathbb{T},
    $$
    then \eqref{localexeqsys} admits global-in-time classical solutions.
    \item
    \emph{Medium alignment} ($\psi_{min}<\lambda\le \psi_{max}$):
    under the admissible condition 
    \begin{equation}\label{eq:medacc}
   \psi_{max} - \psi_{min} < \frac{e^{\frac{\tan^{-1}\hat z}{\hat z}}}{2\left(1 + e^{-\frac{\pi}{\tilde z}} \right)}\lambda,
    \end{equation}
    if the initial data lie in the subcritical region $\Sigma_3$, namely
    $$
   \big(G_0(x),\rho_0(x)\big)\in\Sigma_3, \quad \forall\, x\in\mathbb{T},
    $$
    then \eqref{localexeqsys} admits global-in-time classical solutions.
\end{enumerate}
%Consequently, \eqref{localexeqsys} has a global smooth solution.
Here, the parameters $\hat z$ and $\tilde z$ are defined as
\begin{equation}\label{eq:z}
  \hat z := \sqrt{\left(\frac{\lambda}{\psi_{max}}\right)^2-1}
  \quad\text{and}\quad
  \tilde z := \sqrt{\left(\frac{\lambda}{\psi_{min}}\right)^2-1}.
\end{equation}
Note that $\hat z$, $\tilde z$ could be real, purely imaginary, as
well as infinity.
%The regions $\Sigma_1, \Sigma_2$ and $\Sigma_3$ are subsets of
%$\mathbb{R}\times\mathbb{R}_+$, defined in
%\eqref{invregtransf}, \eqref{invregtransf2} and \eqref{invregtransf3} respectively.
\end{theorem}

The weak and medium alignment situations require an additional structural inequality \eqref{eq:weakacc} and \eqref{eq:medacc} respectively, so that a subcritical region can be obtained through our techniques. These inequalities only depend on the parameters of the EPA system, $k,c,\psi$. These conditions arise due to the presence of oscillatory solutions and as a consequence of our method in handling these to arrive at the thresholds. In the strong alignment case, the solutions decay exponentially to the equilibrium solution without any oscillations, obviating the requirement of any additional condition. 
%It is an open problem as to how the critical thresholds look when \eqref{eq:weakacc} or \eqref{eq:medacc} are not satisfied in the respective cases. Even with this shortcoming, our result signifies a remarkable achievement as threshold for EPA system is quite hard to obtain and our result is the first such result. 

We also prove the corresponding finite-time-breakdown result but do not include it here. We obtain regions $\Delta_1,\Delta_2,\Delta_3$ which are the supercritical regions for the weak, strong and medium alignment cases respectively.

\iffalse
\begin{theorem}[Finite time breakdown]
\label{ftb1}
Under the same assumptions as Theorem \ref{gs1}, we have
\begin{enumerate}
    \item
      \emph{Weak alignment} ($\psi_{max}<\lambda$):
      If there exists $x_0\in\mathbb{T}$ that lie in the supercritical
      region $\Delta_1$, namely 
      $$
      \big(G_0(x_0),\rho_0(x_0)\big)\in\Delta_1,
      $$
      then $(G,\rho)$ becomes unbounded at a finite time.
    \item
      \emph{Strong and medium alignment} ($\psi_{max}\geq\lambda$):
      If there exists $x_0\in\mathbb{T}$ that lie in the supercritical
      region $\Delta_2$, namely 
    $$
        \big(G_0(x_0),\rho_0(x_0)\big)\in\Delta_2,
    $$
    then $(G,\rho)$ becomes unbounded at a finite time.
\end{enumerate}
Moreover, at the blowup time $t_c$ and location $x_c$, the solution
generate a singular shock, with
\[\lim_{t\to t_c^-}  \rho(t,x_c) = \infty \,\,\text{or}\,\,0,\quad
  \lim_{t\to t_c^-} G(t,x_c) = -\infty,\quad
  \lim_{t\to t_c^-}  u_x(t,x_c ) = -\infty.\]
%The regions $\Delta_1, \Delta_2$ are defined in \eqref{ftbinvregtransf}, \eqref{ftbinvregtransf2} respectively.
\end{theorem}
\fi

\begin{figure}[ht!] 
\centering
\subfigure[Weak alignment($\lambda = \sqrt{2}, \psi_{Max} = 0.75,\psi_{min} = 0.25 $)]{\label{subcrGrhoweakinitial} \includegraphics[width=0.45\linewidth]{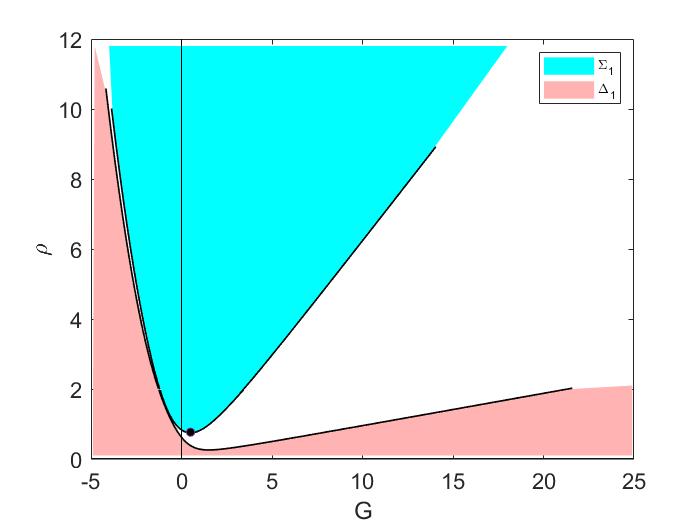}}%
\qquad
\subfigure[Strong alignment ($\lambda = \sqrt{2}, \psi_{max} = 2,\psi_{min} = 1.5 $)]{\label{subcrGrhostronginitial}\includegraphics[width=0.45\linewidth]{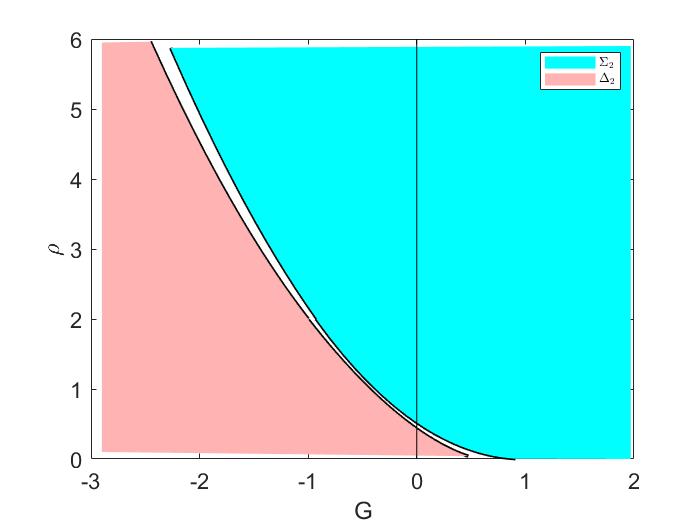}}
\caption{Shapes of $\Sigma_1,\Sigma_2,\Delta_1,\Delta_2$.}
\label{mainfig}
\end{figure}

Similar to the example in the introduction, the fundamental step in analyzing EPA systems for thresholds is to derive an ODE system along the common characteristic path, $\left\{(t,x): \frac{dx}{dt}=u(t,x(t)),\, x(0)=\alpha\right\}$. $\alpha\in\mathbb{T}$ is the parameter which is fixed for a single characteristic path. The resulting ODE system is analyzed for each path. The global well-posedness of unknowns in \eqref{localexeqsys} is obtained by combining the all-time-existence of the unknowns in \eqref{Grhosys} for all $\alpha\in\mathbb{T}$. For \eqref{localexeqsys}, the resulting ODE system is,
\begin{subequations}
\label{Grhosys}
\begin{align}
& G' = -G(G -\psi\ast\rho) + k(\rho -c),  \label{Geq}\\
& \rho' = -\rho (G -\psi\ast\rho). \label{rhoeq}
\end{align}
\end{subequations}

Next, we move on to another important situation of that of the weakly singular kernel, that is, $\psi\in L^1(\mathbb{T})$. Here, we do not have the \eqref{eq:psibounded} type of bounds which are essential in the threshold analysis. 

Therefore, we have to modify our technique. Here, we need to improve the bounds on $\psi\ast\rho$ to obtain valid thresholds. Following is the global existence result.
%The key is to look for plausible constants $\rho_{min},\rho_{max}$ so that $\rho_{min}\leq\rho(\cdot, t)\leq\rho_{max}$ and subsequently find a subcritical region that lies within these bounds on the density, see Figure \ref{figl1}. In this situation, the auxiliary systems to \eqref{mainodesys} are obtained by plugging in improved bounds for $\psi\ast\rho$. The immediate crude bounds for $\psi\ast\rho$, namely $||\psi||_1\rho_{min}$ and $||\psi||_1\rho_{max}$ result in a trivial subcritical region. The analysis is delicate also because the auxiliary system and in turn, the subcritical region themselves depend on the bounds on density, which was not the case for bounded $\psi$. We are able to construct a subcritical region, thereby, proving the existenc of a threshold.

\begin{theorem}[Weakly singular alignment force]  
\label{l1gs}
Let $\psi\in L^1(\mathbb{T})$. Set $\lambda := 2\sqrt{\frac{k}{c}}$. Suppose the initial data to \eqref{localexeqsys} satisfies the hypothesis of Theorem \ref{gs1}. Then there exists sets $\Sigma_L^1,\Sigma_L^2,\Sigma_L^3$ such that,
\begin{enumerate}
    \item \emph{Weak alignment} ($\|\psi\|_{L^1}-\gamma  < \tfrac{\lambda}{2}$): under the admissible condition
    \begin{equation}\label{eq:wscond}
  4(\|\psi\|_{L^1} - 2\gamma) < \frac{e^{\frac{\tan^{-1}\hat z}{\hat z}}\left(1 - e^{ -\frac{\pi}{\tilde z} - \frac{\pi}{\hat z} }\right) }{2\left(1 + e^{-\frac{\pi}{\tilde z}}\right)}\lambda,
  \end{equation}
  if the initial data lie in the subcritical region $\Sigma_L^1$,
namely
\[
 \big(G_0(x),\rho_0(x)\big)\in\Sigma_L^1, \quad \forall\,x\in\mathbb{T},
\]
then $(G,\rho)$ remain bounded in all time.
      \item \emph{Strong alignment} ($\gamma  \geq \tfrac{\lambda}{2}$): if the initial data lie in the subcritical region $\Sigma_L^2$, namely
    $$
   \big(G_0(x),\rho_0(x)\big)\in\Sigma_L^2, \quad \forall\,x\in\mathbb{T},
    $$
    then $(G, \rho)$ remain bounded in all time.
    \item
    \emph{Medium alignment} ($\gamma<\frac{\lambda}{2}\le ||\psi||_1-\gamma$):
    under the admissible condition 
    \begin{equation}\label{eq:wsmedacc}
   4(||\psi||_{L^1}-2\gamma) < \frac{e^{\frac{\tan^{-1}\hat z}{\hat z}}}{2\left(1 + e^{-\frac{\pi}{\tilde z}} \right)}\lambda,
    \end{equation}
    if the initial data lie in the subcritical region $\Sigma_L^3$, namely
    $$
   \big(G_0(x),\rho_0(x)\big)\in\Sigma_L^3, \quad \forall\, x\in\mathbb{T},
    $$
    then $(G, \rho)$ remain bounded in all time.
\end{enumerate}
Consequently, \eqref{localexeqsys} has a global smooth solution.
Here, $\gamma = \int_{1/2}^1\psi^*(x)\, dx$, where $\psi^* :
(0,1]\to\mathbb{R}$ is the decreasing rearrangement of $\psi$ on
$\mathbb{T}$.
The parameters $\hat z$ and $\tilde z$ are defined as
\begin{equation}\label{eq:zl1}
  \hat z := \sqrt{\left(\frac{\lambda}{2(\|\psi\|_{L^1}-\gamma)}\right)^2-1}
  \quad\text{and}\quad
  \tilde z := \sqrt{\left(\frac{\lambda}{2\gamma}\right)^2-1}.
\end{equation}
%The regions $\SigmaL^1,\SigmaL^2$ and $\SigmaL^3$ are subsets of $\mathbb{R}\times\mathbb{R}_+$ defined in \eqref{l1subcr}, \eqref{l1subcr2} and \eqref{l1subcr3} respectively.
\end{theorem}

\begin{figure}[ht!] 
\centering
\subfigure[Weak alignment ($\lambda = 4, \|\psi\|_{L^1} = 2, \gamma=0.95$)]{\label{l1subcrfig} \includegraphics[width=0.45\linewidth]{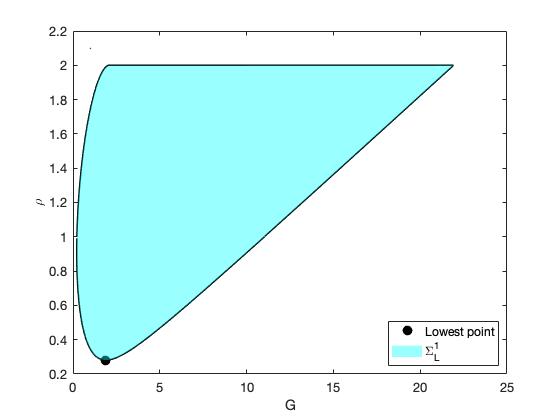}}%
\qquad
\subfigure[Strong alignment ($\lambda = \sqrt{2}, \|\psi\|_{L^1} = 2, \gamma=0.95 $)] {\label{wssalsubcr}\includegraphics[width=0.45\linewidth]{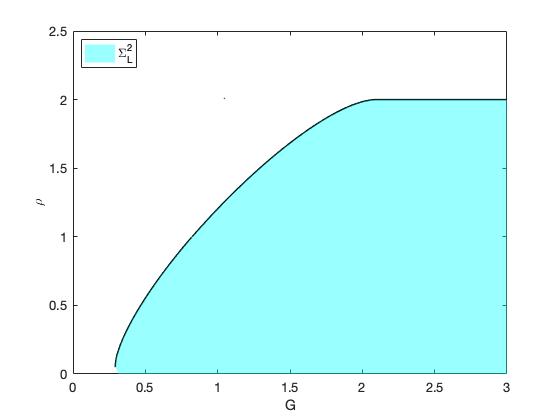}}
\caption{Shapes of $\Sigma_L^1,\Sigma_L^2$.}
\label{figl1}
\end{figure}

%Unlike the case when $\psi$ is bounded, the subcritical regions
 % $\Sigma_L^i$'s here are subsets of $\{(G_0,\rho_0) :
%  \rho_{min}\le\rho_0\le\rho_{max}\}$ for appropriate choices of
 % $0\le\rho_{min}<c<\rho_{max}<\infty$.
  Figure \ref{figl1} illustrates the
  shape of $\Sigma_L^1$ and $\Sigma_L^2$.
  The steady-state solution $(G,\rho)=(c\|\psi\|_{L^1},c)\in\Sigma_L^i$.
  Hence, the region $\Sigma_L^i$ contains initial data around the steady state and we obtain a nontrivial subcritical region.

  The admissible conditions \eqref{eq:wscond} and \eqref{eq:wsmedacc} are similar to
  \eqref{eq:weakacc} and \eqref{eq:medacc} respectively, in the sense that both imply that nonlocality of $\psi$ is not too strong. When $\psi$ is bounded, $\psi_{max}-\psi_{min}$ is its oscillation. It is zero if and only if $\psi$ is constant. Correspondingly, when $\psi$ is unbounded but integrable, $\psi_{max}-\psi_{min}$ is
  replaced by $4(\|\psi\|_{L^1}-2\gamma)$. Note that
  $\|\psi\|_{L^1}-2\gamma\ge 0$, and the equality holds if and only if
  $\psi$ is a constant. Therefore, all of \eqref{eq:wscond}, \eqref{eq:wsmedacc}, \eqref{eq:weakacc}, \eqref{eq:medacc} impose an upper bound on how much $\psi$ is offset from a constant.
  %Hence, just like in the $L^\infty$ case, the admissible condition says that the
%The parameters $\hat z$ and $\tilde z$ are also revised to adapt the unboundedness of $\psi$.

The following result is for the EA system with weakly singular influence function.
\begin{theorem}
\label{eal1}
Let $\psi\in L^1(X)$ ($X=\mathbb{R}$ or $\mathbb{T}$) be non-negative. Consider \eqref{localexeqsys} with $k=0$. If 
$$
\inf_x G(0, x)>0,
$$
then there exists global-in-time $C^1$ solution to the system \eqref{localexeqsys}. Moreover, $\rho,G$ have uniform bounds in terms of $\psi,\rho_0,G_0$. 
\end{theorem}

\section{Nonlocal Euler system with relaxation}
\label{sec:3}
Several traffic-flow and fluid-flow models are modeled through a density that follows the continuity equation with a nonlocal flux, see for example \cite{FV16,LT22}. A general equation is,
$$
\rho_t + (\rho \text{v})_x = 0,
$$
for some nonlocal $\text{v}$. We augment this model with a velocity obeying Burgers' equation with some source terms. The source terms are such that the particles modeled move towards an equilibrium state, imparting some order and regularity to the system. We consider the following pressureless Euler-like model,
\begin{subequations}
\label{hypMain}
\begin{align}
& \rho_t +(\rho \text{v})_x=0,\; x\in \mathbb{R}, \; t>0, \label{hypMainmass}\\
&  u_t+ uu_x=\rho(\text{v}-u), \label{hypMainbur}
\end{align}
\end{subequations}
with initial data $(\rho_0\geq 0, u_0)$. 
%To close this system, we need to relate $\text{v}$ to the unknown variables, $\rho,u$. 
The quantity $\text{v}$ determines the steady state velocity. 
%Such a system can be considered as the augmentation of the continuity equation with an equation for the velocity of particles. 
Motivated by physical assumptions, we set $\text{v}=Q\ast u$, with,
\begin{align}
\label{qhyp}
\begin{aligned}
& Q\in W^{1,1}(\mathbb{R}),\quad Q(x) = Q(-x) \text{ (symmetric)}, \\
& \int_{\mathbb{R}}\!\! Q(x)\, dx = 1,\quad \text{and } Q \text{ is decreasing away from origin}.
\end{aligned}
\end{align}
Evidently, $Q\ast u$ is a weighted average of the velocity, with maximum weight at the point itself and decreasing symetrically as one moves away. To our knowledge, the well-posedness of \eqref{hypMain} was not known. Due to the nonlocal flux, it cannot be concluded from existing literature for hyperbolic balance laws. Inspired by Kato and Majda, see \cite{Ka75,Ma84}, we use energy methods to prove local existence/uniqueness for a general multidimensional system, the one dimensional form of which is \eqref{hypMain}, in a relatively more general space. In particular, we allow for solutions that need not decay at infinity, though they are bounded. We do not state the result here but instead, focus on the one dimensional threshold result.
\begin{theorem}
\label{finalresult}
Consider \eqref{hypMain} with $\text{v}=Q\ast u$ and $Q$ satisfying \eqref{qhyp}. Suppose the initial density is nonnegative and the initial data $\rho_0,u_0\in L^\infty(\mathbb{R})$, and $\rho_{0x},u_{0x}\in H^s(\mathbb{R})$, $s\geq 1$.
\begin{itemize}
\item{[Subcritical region]} A unique solution $\rho,u\in C([0,\infty);L^\infty(\mathbb{R}))$ and 
$$\rho_x , u_x \in C([0,\infty);H^s(\mathbb{R})), \quad s\geq 1 $$ exists if $u_{0x}(x) + \rho_0 (x) \geq 0$ for all $x\in\mathbb{R}$. 
\item{[Supercritical region]} If $\exists x_\ast\in\mathbb{R}$ for which $u_{0x}(x_\ast) < -\rho_0 (x_\ast)$, then $u_x\to-\infty$ in finite time.
\end{itemize}
\end{theorem}

Another interesting case is when $\text{v}$ is local in \eqref{hypMain}, that is, $\text{v}=f(\rho,u)$. Let $\Theta = \{(\rho,u): u=f(\rho,u)\}$. Then \eqref{hypMain} is strictly hyperbolic as long as $(\rho,u)\in\Theta^c$. However, this cannot be guaranteed a priori and \eqref{hypMain} could degenerate from strict hyperbolic to weak hyperbolic at certain time. It turns out that for $f_\rho = 0$ ($f$ only depending on the velocity), strict hyperbolicity can be guaranteed if the initial data lies in a certain set. The threshold analysis and results depend heavily on whether \eqref{hypMain} is strictly or weakly hyperbolic.

\begin{theorem}
\label{thm1}
Consider the system \eqref{hypMain} with $\text{v}=f(u) $ and initial conditions $ (\rho_0\geq 0 ,u_0) \in C_b^1(\mathbb{R})\times C_b^1(\mathbb{R})$ with $ \inf|f(u_0) - u_0 |>0$.  If $f_u \leq 0$ for solution $u$ under consideration, then, 
\begin{enumerate}
\item
\textbf{Bounds on $u$ and $\rho$:}
$u(t,\cdot)$ is uniformly bounded and satisfies $|f(u(t,\cdot))-u(t,\cdot)|>0$ for $t>0$. And
$$
\rho(t,x)\leq \frac{\sup\rho_0|f(u_0)-u_0|}{e^{\int_{u_0(x)}^{u(t,x)}\frac{d\xi}{f(\xi)-\xi}}|f(u(t,x))-u(t,x)|}.
$$
\item
\textbf{Global solution:} If
$$
u_{0x}(x)+\rho_0(x)\geq 0,\quad \forall x\in\mathbb{R},
$$ 
then there exists a global classical solution $\rho,u\in C^1((0,\infty)\times\mathbb{R})$.
Moreover, $\rho, u_x$ are uniformly bounded  with 
$$
0 \leq \rho(t,x) \leq M, \quad 0 \leq u_x(t,x) + \rho(t,x)  \leq M, \qquad \forall t>0,x\in\mathbb{R},
$$
where $M = \max\{ \sup \rho_0 ,\sup (u_{0x}+\rho_0) \}$. 
\item
\textbf{Finite time breakdown:} If $\exists x_0\in\mathbb{R}$ such that
$$
u_{0x}(x_0)+ \rho_0(x_0)<0,
$$
then $\lim_{t\to t_c} |\rho_x| =\infty$ or $\lim_{t\to t_c} |u_x| = \infty$ for some $t_c>0$.
\end{enumerate}
\end{theorem}

The condition $\inf|f(u_0) - u_0 |>0$ ensures strict hyperbolicity by ensuring $\inf|f(u(t,\cdot)) - u(t,\cdot) |>0$ for all $t>0$. If this does not hold, then we have Theorem \ref{thm2}. Also, the upper bound on $\rho$ can be infinitely large as the system gets `closer' to weakly hyperbolic. This indeed exhibits the borderline behaviour of density between strictly and weakly hyperbolic systems. In strictly hyperbolic systems (with non-erratic source terms) density is bounded for all times, even when shock forms, whereas in weakly hyperbolic systems, density becomes unbounded when shock forms, as is the case in \eqref{epamain}.

\begin{theorem}
\label{thm2}
Let $f $ be a smooth function depending on $u$ only, i.e., $f_\rho = 0$. Consider the system \eqref{hypMain} subject to initial conditions, $
(\rho_0\geq 0 ,u_0) \in C^1_b(\mathbb{R})\times C^2_b(\mathbb{R}) $.
%$s\geq 1$ is an integer. 
 If $f_u \leq 0$ for the solution $u$ of consideration, then 
\begin{enumerate}
%\item
%\textbf{Bounds on $u$:}
%$u(\cdot,\cdot)$ is uniformly bounded. 
\item
\textbf{Global Solution:} If
$$
u_{0x}(x)+\rho_0(x)\geq 0,\quad \forall x\in\mathbb{R},
$$
then there exists a global solution 
$$
\rho, u\in C^1((0,\infty)\times \mathbb{R}).
$$  
Moreover, $u,\rho, u_x$ are uniformly bounded.
%and,
%$$
%0 \leq \rho(t,x) \leq M, \quad 0 \leq u_x(t,x) + \rho(t,x) \leq M, \qquad \forall t>0,x\in\mathbb{R},
%$$
%where $M = \max\{ \sup \rho_0 ,\sup (u_{0x}+\rho_0) \}$. And 
Also, we have the following,
$$
||\rho_x(t,\cdot)||_\infty \leq C_1 e^{C_2t},\quad t>0,
$$
where $C_1 = C_1(||\rho_0||_{C^1}, ||u||_{C^2})$ and $C_2 = C_2(||u||_{C^2},||\rho_0||_\infty)$.
\iffalse
\item
\textbf{Finite time breakdown:} If $\exists x_0\in\mathbb{R}$ such that
$$
u_{0x}(x_0)+ \rho_0(x_0)<0,
$$
then $\lim_{t\to t_c} u_x = -\infty$ for some $t_c>0$ at the rate of $O\left(\frac{1}{|t-t_c|}\right)$ or faster.
\fi
\end{enumerate}
\end{theorem}

We would like to point out some key differences in Theorems \ref{thm1} and \ref{thm2}. Firstly, $\rho,u$ are bounded for all times in the former which is not true for the latter wherein there might be density concentration, that is, $\rho\to\infty$ in finite time. Secondly, the space in which the solutions lie is different due to an extra degree of smoothness needed for velocity which arises in proving the local existence in the case of pressureless Eulerian systems.

The following result is for a general $f$, that is, it is a function of both density and velocity. Here, we cannot guarantee a priori strict hyperbolicity and it is imperative that we consider \eqref{hypMain} as weakly hyperbolic. As a result, we need more conditions for global existence.
\begin{theorem}
\label{thm3}
Let $f = f(\rho,u) $ be a smooth function of its variables. Consider the system \eqref{hypMain} with initial conditions $(\rho_0\geq 0 ,u_0) \in C_b^1(\mathbb{R})\times C_b^2(\mathbb{R})$. If $f_u \leq 0$ for the solutions under consideration, then $u,\rho,u_x$ are uniformly bounded. 
\iffalse
\begin{enumerate}
\item
\textbf{Bounds on $u$:} There exists a smooth function $\phi:\mathbb{R}^+\to\mathbb{R}$ such that,
$$
\min\left\{ \inf_{\mathbb{R}} u_0, \min \phi(\rho) \right\}\leq u(t,\cdot )\leq \max\left\{ \sup_{\mathbb{R}} u_0, \max\phi(\rho) \right\},
$$
for as long as $\rho \geq 0$ is bounded.
\item
\textbf{Bounds on $\rho,u_x$:} If
$$
u_{0x}(x)+\rho_0(x)\geq 0,\quad \forall x\in\mathbb{R},
$$
then $\rho, u_x$ are uniformly bounded  with 
$$
0 \leq \rho(t,x) \leq M, \quad 0 \leq u_x(t,x)+ \rho(t,x)   \leq M, \qquad \forall t>0,x\in\mathbb{R}.
$$
where $M = \max\{ \sup \rho_0 ,\sup (u_{0x}+\rho_0) \}$. 
\item
\fi
If in addition to $u_{0x}+\rho_0\geq 0$,
\begin{itemize}
    \item $(\rho f)_{\rho\rho}\geq 0$, $f_{uu}\leq 0$ along with 
    $$
    \rho_{0x}(x)\geq 0,\quad  u_{0xx}(x)+\rho_{0x}(x)\geq 0,\ \forall x\in\mathbb{R},
    $$ 
    \begin{center} OR \end{center}
    \item $(\rho f)_{\rho\rho}\leq 0$, $f_{uu}\geq 0$ along with 
    $$
    \rho_{0x}(x)\leq 0,\quad u_{0xx}(x)+\rho_{0x}(x)\leq  0,\ \forall x\in\mathbb{R},
    $$
\end{itemize}
then there exists a global solution $\rho, u\in C^1((0,\infty)\times\mathbb{R})$. In addition,
$$
u_x + \rho \in C^1((0,\infty)\times\mathbb{R}).
$$
\iffalse
\item
\textbf{Finite time breakdown:} If $\exists x_0\in\mathbb{R}$ such that
$$
u_{0x}(x_0)+ \rho_0(x_0)<0,
$$
then $\lim_{t\to t_c} u_x = -\infty$ at the rate of $O\left(\frac{1}{|t-t_c|}\right)$ or faster, for some $t_c>0$.
\end{enumerate}
\fi
\end{theorem}

The breakdown result is same for Theorems \ref{thm2} and \ref{thm3}.
\begin{theorem}[Finite time breakdown]
\label{thm4}
Consider \eqref{hypMain} with $\text{v}=f(\rho,u)$. If there exists an $x_\ast$ such that $u_{0x}(x_\ast)+\rho_0(x_\ast) < 0$, then $\lim_{t\to t_c^-}|u_x(t,x_c)| = -\infty$ for some $x_c$ and $t_c>0$. 
\end{theorem}

%Owing to the source term in \eqref{hypMainbur}, which ensures that particles are pulled towards the equilibrium state, we can derive a maximum principle-like result for the velocity. Therefore, uniform bounds on $u$ are essentially a direct result of the equation system.

A key step is to identify a quantity $e:=u_x+\rho$, that simplifies \eqref{hypMain}. $e$ is quite analogous to $G$ in \eqref{localexeqsys}. This transformation results in,
\begin{subequations}
\label{coupleerho}
\begin{align}
& \rho_t + f\rho_x = f_u\rho(\rho - e), \label{coupleerho1}\\
& e_t + u e_x = -e(e-\rho). \label{coupleerho2}
\end{align}	
\end{subequations}
If $f_u\leq 0$, we can bound $e,\rho$ in tandem, that is, both quantities are all-time-bounded or break down together. Theorems \ref{thm2} and \ref{thm3} can be proved thereafter.

Under the assumptions of Theorem \ref{thm1}, where the system is strictly hyperbolic, we find the two Riemann invariants and analyze them. From \eqref{hypMain}, one of the Riemann invariants is simply $u$. The other invariant $R$ can be evaluated using conventional techniques.
%and is,
%\begin{align*}
%\label{R1}
%R(\rho,u) := \rho e^{\int_{u_0}^u \frac{d\xi}{f(\xi)-\xi}} (f(u)-u).
%\end{align*}
Bounds on $R_x$ ensure bounds on $\rho_x$. The asymptotic bound on $\rho$ in Theorem \ref{thm1} is obtained by bounding $R$.

\section*{Acknowledgments}
This work was supported in part by the National Science Foundation under
Grant DMS1812666.

\bigskip

\bibliographystyle{abbrv}

\end{document}